\documentclass[a4paper]{jpconf}
\usepackage{graphicx}

\usepackage{amsfonts,amsmath,amscd,amssymb,amsthm}
\def\0{{\bar 0}}
\def\1{{\bar 1}}

\def\Z{{\mathbb Z}}

\def\N{{\mathbb N}}

\def\hgt{{\operatorname{ht}}}

\def\ch{{\operatorname{ch}\:}}

\newcommand{\ttk}{\mathtt{k}}

\newcommand{\gL}{\Lambda}

\newcommand{\itema}{\indent {{\rm(a)}}}
\newcommand{\itemb}{\indent {{\rm(b)}}}
\newcommand{\itemi}{\indent{{\rm(i)}}}

\newcommand{\itemc}{\indent{{\rm(c)}}}
\newcommand{\itemd}{\indent{{\rm(d)}}}
\newcommand{\iteme}{\indent{{\rm(e)}}}
\newcommand{\itemf}{\indent{{\rm(f)}}}
\newcommand{\itemg}{\indent{{\rm(g)}}}
\newcommand{\itemh}{\indent{{\rm(h)}}}

\newcommand{\noi}{\noindent}
\newcommand{\ga}{\alpha}
\newcommand{\gb}{\beta}
\newcommand{\gc}{\gamma}

\newcommand{\gd}{\delta}

\newcommand{\gs}{\sigma}

\newcommand{\gt}{\tau}

\newcommand{\gl}{\lambda}

\newcommand{\gr}{\rho}

\newcommand{\gep}{\epsilon}
\newcommand{\gth}{\theta}

\newcommand{\ot}{\otimes}

\newcommand{\fg}{\mathfrak{g}}
\newcommand{\fsl}{\mathfrak{sl}}\newcommand{\osp}{\mathfrak{osp}}

\newcommand{\fh}{\mathfrak{h}}

\newcommand{\fb}{\mathfrak{b}}

\newcommand{\fn}{\mathfrak{n}}

\newcommand{\ff}{\footnote}
\newfont{\eufm}{eufm10 scaled\magstep1}

\newcommand{\cO}{\mathcal{O}}

\newcommand{\cH}{\mathcal{H}}

\newcommand{\ey}{\end{eqnarray}}
\newcommand{\by}{\begin{eqnarray}}

\newcommand{\bco}{\begin{conjecture}}
\newcommand{\ba}{\begin{alg}}
\newcommand{\ea}{\end{alg}}
\newcommand{\eco}{\end{conjecture}}
\newcommand{\bpf}{\begin{proof}}
\newcommand{\epf}{\end{proof}}
\newcommand{\bt}{\begin{theorem}}

\newcommand{\et}{\end{theorem}}

\newcommand{\er}{\end{rem}}
\newcommand{\brs}{\begin{rems}}
\newcommand{\ers}{\end{rems}}
\newcommand{\bl}{\begin{lemma}}
\newcommand{\bsul}{\begin{sublemma}}
\newcommand{\esul}{\end{sublemma}}
\newcommand{\bp}{\begin{proposition}}
\newcommand{\be}{\begin{equation}}
\newcommand{\bc}{\begin{corollary}}
\newcommand{\bexs}{\begin{examples}}
\newcommand{\eexs}{\end{examples}}
\newcommand{\bexa}{\begin{example}}
\newcommand{\eexa}{\end{example}}
\newcommand{\bex}{\begin{exercise}}
\newcommand{\eex}{\end{exercise}}
\newcommand{\btab}{\begin{tab}}
\newcommand{\etab}{\end{tab}}
\newcommand{\el}{\end{lemma}}
\newcommand{\ep}{\end{proposition}}
\newcommand{\ee}{\end{equation}}
\newcommand{\ec}{\end{corollary}}
\newcommand{\Bc}{\begin{center}}
\newcommand{\Ec}{\end{center}}
\newcommand{\bh}{\begin{hyp}}
\newcommand{\eh}{\end{hyp}}
\newcommand{\bhs}{\begin{hyps}}
\newcommand{\ehs}{\end{hyps}}
\newcommand{\bd}{\begin{dfn}}
\newcommand{\ed}{\end{dfn}}

\newtheorem{thm}{Theorem}[section]
\newtheorem{hyp}[thm]{Hypothesis}
 \newtheorem{hyps}[thm]{Hypotheses}
  \newtheorem{rems}[thm]{Remarks}
\newtheorem{conjecture}[thm]{Conjecture}
\newtheorem{theorem}[thm]{Theorem}
\newtheorem{theorem a}[thm]{Theorem A}
\newtheorem{example}[thm]{Example}
\newtheorem{examples}[thm]{Examples}
\newtheorem{corollary}[thm]{Corollary}
\newtheorem{rem}[thm]{Remark}
\newtheorem{lemma}[thm]{Lemma}

\newtheorem{proposition}[thm]{Proposition}

\newtheorem{exercise}[thm]{Exercise}
\newcommand{\lra}{\longrightarrow}

\begin{document}

\title{\v Sapovalov elements and the Jantzen filtration for contragredient Lie superalgebras: A Survey.}

\author{Ian M. Musson}

\address{Department of Mathematical Sciences,
University of Wisconsin-Milwaukee, USA.}

\ead{musson@uwm.edu}

\begin{abstract}
This is a survey of some recent results on
\v Sapovalov elements and the Jantzen filtration for contragredient Lie superalgebras.
The topics covered include the  existence and uniqueness of the  \v Sapovalov elements, bounds on the degrees of their coefficients and the
 behavior of \v Sapovalov elements
  when the Borel subalgebra is changed.  There is always a unique term whose coefficient has larger degree than any other term. This allows us to define some new highest weight modules. If $X$ is a set of orthogonal isotropic roots and $\gl \in \fh^*$ is such that $\gl +\gr$ is orthogonal to all
roots in $X$, 	we construct highest weight modules with character $\gep^\gl{p}_X$. Here $p_X$ is a partition function that counts partitions not involving roots in $X$. When $|X|=1,$  these modules are used to give a  Jantzen sum formula for Verma modules in which all terms are characters of modules in the category $\cO$ with positive coefficients.

\end{abstract}

\section{Introduction.} \label{u.1}

Throughout this paper we work over an  algebraically closed field $\ttk $
  of characteristic zero. If  $\fg$ is a semisimple Lie algebra,
necessary and sufficient conditions for the existence of a non-zero homomorphism between Verma modules can be obtained by combining work of Verma \cite{Ve} with work of Bernstein, Gelfand and  Gelfand \cite{BGG1}, \cite{BGG2}.  Such maps can be described explicitly in terms of certain elements introduced by   \v Sapovalov in \cite{Sh}. Necessary and sufficient  conditions for a simple highest weight module to be a composition factor of a Verma module were also obtained in \cite{BGG1}, \cite{BGG2}. A more elementary proof of the latter result can be given using the Jantzen filtration and sum formula \cite{J1}. \ff{ We remark that the multiplicities of the composition factors in Verma modules for semisimple Lie algebras are given by the Kazdhan-Lusztig conjecture \cite{KaLu}. For type A Lie superalgebras, they are given by Brundan's analog of the Kazdhan-Lusztig conjecture, proven in \cite{CLW}, see also \cite{BLW} Theorem A and \cite{B} Theorem 3.6.}
However neither \v Sapovalov elements
nor the  Jantzen filtration have received much attention
for classical simple Lie superalgebras.
In  this paper we review some recent results on
\v Sapovalov elements, the Jantzen filtration  and sum formula in the super case. New phenomena arise due to the presence of isotropic roots. Proofs not found here are given in \cite{M} Chapter 9 and 10, \cite{M1} or \cite{M17}.
\\ \\
 Let $\fg= \fg(A,\gt)$ be a finite dimensional contragredient Lie superalgebra with Cartan subalgebra $\fh$, and set of simple roots
$\Pi$, see \cite{M} Chapter 5.
	Let $\Delta^{+}$ be the set of positive roots  containing $\Pi$, and  %

%\be \label{gtri}\mathfrak{g} = \mathfrak{n}^- \oplus \mathfrak{h} \oplus \mathfrak{n}^+\ee
\[\mathfrak{g} = \mathfrak{n}^- \oplus \mathfrak{h} \oplus \mathfrak{n}^+\]
the corresponding triangular decomposition  of $\fg$. We use the Borel subalgebras $\mathfrak{b} =  \mathfrak{h}
\oplus \mathfrak{n}^+$ and
$\fb^- =
\mathfrak{n}^- \oplus \mathfrak{h}$.
 The Verma module
$M(\gl)$ with highest weight $\gl \in \fh^*$, and highest weight vector
$v_\lambda$
is induced from $\mathfrak{b}$.
\noi Let $\gr_0$ (resp.
$\gr_1$) be the half-sum of the positive even (resp. odd) positive roots and $\gr=
\gr_0 -\gr_1.$
\\ \\
\noi Fix a non-degenerate invariant symmetric bilinear form $(\;,\;)$ on $\fh^*$ as in \cite{M} Theorem 5.4.1, and for all $\ga \in \fh^*$, let $h_\ga \in \fh$ be the unique element such that $(\ga,\gb) = \gb(h_\ga)$ for all $\gb \in \fh^*$.
Then for all $\alpha
\in \Delta^+$, choose elements $e_{\pm \alpha} \in
\mathfrak{g}^{\pm \alpha}$
 such that

\[ [e_{\alpha}, e_{-\alpha}] = h_{\alpha}.\]
%\be \label{ebl} e_\ga e_{-\ga} v_\mu = h_\ga v_\mu =(\mu, \ga)v_\mu . \ee
  We give bounds on the degrees of the coefficients
of $H_{\pi}$ of the \v Sapovalov elements  in Equation (\ref{rat}) below.  The exact form of the coefficients depends on the way the positive roots are ordered. However there is always a unique coefficient  of highest degree, and we determine the leading term of this coefficient up to a scalar multiple. These results appear to be new even for simple Lie algebras.
 The existence of a unique coefficient  of highest degree is used to construct some new highest weight modules $M^\gc(\gl)$, where $\gc$ is an isotropic root and $(\gl+\gr,\gc) =0$, \cite{M17}.  There is a version of the Jantzen sum formula in \cite{M} Theorem 10.3.1, see also \cite{G},  which involves infinite sums of Verma modules with alternating sign coefficients.  The module  $M^\gc(\gl)$ has character  $\epsilon^{\lambda}p_\gc$ see (\ref{sing}) for notation. The infinite sums can be replaced by   a finite sum  of characters of modules of this form leading to a formula where both sides are sums of characters in the category $\cO.$ In type A there are explicit expressions for \v Sapovalov elements, see \cite{M1} Section 7.
%I  would like to thank %Jon Brundan for suggesting the use of noncommutative determinants to write \v Sapovalov elements in Section \ref{s.8}, and raising the possibility of using Theorem \ref{shgl} to prove Theorem \ref{shapel}. I alsothank Kevin Coulembier for some helpful conversations.

\section{\bf Preliminaries.}  \label{sss7.1}
{\bf To simplify the exposition, we assume in this survey that any non-isotropic root of $\fg$ is even.}  This should give the reader the flavor of the results in general, but the proofs and some of the statements
should be modified if this is not the case. \\ \\
If $\eta \in Q^+=\sum_{\ga\in \Pi}\mathbb{N}\ga$, a {\it
partition} of $\eta$ is a map
$\pi: \Delta^+ \longrightarrow
\mathbb{N} $ such that
$\pi(\alpha) = 0$ or $1$ for all isotropic roots $\alpha$,
 and
\[ \sum_{\alpha \in {{\Delta^+}}} \pi(\alpha)\alpha = \eta.\]

\noi Let $\bf{ P(\eta)}$ be the set of partitions of $\eta$.
The {\it degree} of a partition $\pi$ is defined to be $|\pi| = \sum_{\alpha \in \Delta^+} \pi(\alpha).$
Partitions are useful because they can be used to index a basis for $U(\mathfrak{n}^-)$.
 Fix an ordering on the set $\Delta^+$, and for $\pi$ a partition,
set

%\be \label{negpar} e_{-\pi} = \prod_{\alpha \in \Delta^+} e^{\pi (\alpha)}_{-\alpha},\ee
\[ e_{-\pi} = \prod_{\alpha \in \Delta^+} e^{\pi (\alpha)}_{-\alpha},\]
the product being taken with respect to this order.
Then the elements $e_{- \pi},$ with
$\pi \in \bf{ P}(\eta)$ form a basis of
$U(\mathfrak{n}^-)^{- \eta}.$
\\ \\Suppose $\gc$ is a positive root, $m$ is a  positive   integer and set

\[{\mathcal H}_{\gc, m} = \{ \lambda \in  {\mathfrak h}^*|(\lambda + \rho, \gc) = m(\gc, \gc)/2  \},\]
and let ${\mathcal I}(\mathcal H_{\gc, m})$ be the ideal of $S(\fh)$ consisting of  functions vanishing on $\mathcal{H}_{\gamma,m}$. The {\it \v Sapovalov element}
$\theta_{\gamma,m}$
 corresponding to the pair $(\gc, m)$ has the form

\be \label{rat}\theta_{\gamma,m} = \sum_{\pi \in {{\bf P}}(m\gamma)} e_{-\pi}
H_{\pi},\ee
where $H_{\pi} \in U({\mathfrak h})$, and
satisfies

\be \label{boo} e_{\ga} \theta_{\gamma,m} \in U({\mathfrak g})(h_{\gc} + \rho(h_{\gc})-m(\gc,\gc)/2)+U({\mathfrak g}){\mathfrak n}^+ , \; \rm{for \; all }\;\ga \in \Delta^+. \ee
This means that if $\gl\in \cH_{\gc, m}$, then
$\theta_{\gc, m}v_{\gl}$ is a highest weight vector  in $M(\gl)$.
If $\gc$ is simple, then $\theta_{\gamma,m}=e^m_{-\gamma}$, so now assume that $\gc$ is not simple. Let
$\pi^0 \in {{\bf P}}(m\gc)$ be the unique partition of $m\gc$ such
that $\pi^0(\ga) = 0$ if $\ga \in \Delta^+ \backslash \Pi.$
The partition $m\pi^{\gamma}$ of $m\gc$ is given by $m\pi^{\gamma}(\gc)=m,$ and $m\pi^{\gamma}(\ga)= 0$ for all positive roots $\ga$ different from $\gc.$
We normalize  $\theta_{\gamma,m}$ so that the coefficient $H_{\pi^0}$ is equal to 1.
This guarantees that $\theta_{\gc, m}v_\gl$ is never zero.
For a semisimple Lie algebra, the existence of such elements was shown by \v Sapovalov, \cite{Sh} Lemma 1.
\\ \\
\noi For a non-isotropic root $\ga,$ we set $\alpha^\vee = 2\alpha /
(\alpha, \alpha)$, and write $
s_\ga$ for the reflection corresponding to $\alpha$.
As usual the  Weyl group  $W$ is the subgroup of $GL({\fh}^*)$
generated by all such reflections.  For $u \in W$ set

 $$ N(u) = \{ \alpha \in \Delta_0^+ | u \alpha < 0 \},\qquad \ell(u) = |N(u)|.$$
 Suppose $\Pi_{\rm even}$ is the set of  even simple roots, and let $W_{\rm even}$ be
the subgroup of $W$ generated by the reflections $s_\ga,$ where $\ga
\in \Pi_{\rm even}$.
Note that $W_{\rm even}$ can be a proper subgroup of $W$.  For example this happens when $\fg = \osp(2m,2n).$
If $\Pi = \{\ga_i|i   = 1, \ldots, t \}$ is the set of simple roots, and $\gc$ is a positive root such that $\gc
= \sum_{i=1}^t a_i\ga_i,$ then the {\it height} $\hgt \gc$ of $\gc$ is
$\hgt \gc = \sum_{i=1}^t a_i$.
\begin{theorem} \label{Shap}
Let  $\fg$ be semisimple Lie algebra or a contragredient Lie superalgebra, and $\gc$ a positive root.  If $\gc$ is isotropic assume that $m=1.$
Suppose $\gc = w\gb$ \mbox{ for a simple root } $\gb$ \mbox{ and } $w \in W_{\rm even},$ and for $\alpha \in N(w^{-1}),$ set $q(w,\ga) = (w\gb, \ga^\vee).$
Then \\ \\
\itema  \; there exists a \v Sapovalov  element
$\theta_{\gamma,m} \in U({\mathfrak b}^{-} )^{- m\gamma}$, which is unique modulo the left ideal
$U({\mathfrak b}^{-} ){\mathcal I}(\mathcal{H}_{\gamma,m})$.
\\
\\
\itemb \;
the coefficients of $\theta_{\gc, m}$ satisfy

\be \label{x1} |\pi|+\deg  H_{\pi} \le m\hgt \gc, \ee and
%\[ \label{x1} |\pi|+\deg  H_{\pi} \le m\hgt \gc, \] and
%\[\label{hig} H_{m\pi^{\gamma}} \mbox{has leading term } \prod_{\ga \in N(w^{-1})}h_\ga^{mq(w,\ga)}.\]

\be\label{hig} H_{m\pi^{\gamma}} \mbox{has leading term } \prod_{\ga \in N(w^{-1})}h_\ga^{mq(w,\ga)}.\ee
\end{theorem}
\noi  Let $d_{m\gc}$ be the degree of $H_{m\pi^{\gamma}}$.
\bc \label{zoo} In Theorem \ref{Shap} $H_{m
\pi^{\gamma}}$ is the unique term of degree $d_{m\gc}$ in $\gth_{\gc, m}.$
\ec
\bpf This follows easily from \eqref{x1} and \eqref{hig}.\epf
\section{\v Sapovalov Elements and their Coefficients.} \label{s.2}

 \noi A finite dimensional contragredient Lie superalgebra $\fg$ has, in general several conjugacy classes of Borel subalgebras, and this both complicates and enriches the representation theory of $\fg$.  The complications are partially resolved by at first fixing a Borel subalgebra $\fb$ (or equivalently a basis of simple roots for $\fg$) with special properties. The effect of changing the Borel subalgebra is studied in detail in \cite{M17}, see also the next Section.
\\ \\
In \cite{K} Table VI, Kac gave a particular diagram in each case that we will call {\it distinguished.} The corresponding set of simple roots and Borel subalgebra are also called distinguished. The distinguished Borel subalgebra contains at most one simple isotropic root vector.  Unless $\fg = \osp(1,2n)$, or $\fg = \osp(2,2n)$ there is exactly one other Borel subalgebra with this property up to conjugacy in Aut $\fg$. A representative of this class (and its set of simple roots) will be called {\it anti-distinguished.}
We assume that $\fb$ is either  distinguished  or anti-distinguished.
\\ \\
\noi Theorem \ref{Shap} is  proved by looking at the proofs given in \cite{H2} or \cite{M} and keeping track of the coefficients. %
Given $\gl \in {\mathfrak h}^* $ and $x =
\sum_{i} a_i \otimes b_i \in U(\fb^-)= U({\mathfrak n}^-)\otimes S({\mathfrak h}),$ set $x(\gl) =\sum_{i} a_i  b_i(\gl)\in U(\fn^-).$
Let  $(\gc,m)$ be as in the statement of the Theorem and set ${\mathcal H} = {\mathcal H}_{\gc, m}$.
The idea of the proof is to
construct elements $\theta^\gl \in U({\mathfrak n^-})^{-m\gamma}$
for all $\gl$ in a dense subset of ${\mathcal H} $
such that $\theta^\gl v_\gl$ is  a highest weight vector in
$M(\lambda)^{\gl-m\gamma}$, and that

\be \label{yam}\theta^\gl = \sum_{\pi \in {{\bf P}}(m\gc)}a_{\pi, \gl}e_{-\pi}.\ee
where $a_{\pi, \gl} $ is a polynomial function of $\gl \in \Lambda$ satisfying suitable conditions.
For $\pi \in {{\bf P}}(m\gc)$,
the assignment $\gl
 \rightarrow a_{\pi, \gl}$ for $\gl \in \Lambda$ determines a polynomial map from ${\mathcal H}$ to
 $U({\mathfrak n}^-)^{-m\gamma},$ so  there exists  an element  $H_\pi \in U({\mathfrak h})$ uniquely determined modulo ${\mathcal I}(\mathcal H)$
such that $H_{\pi}(\gl) = a_{\pi, \gl}$ for all $\gl \in \Lambda$.   We define the element
$\theta \in U(\fb^-)$  by
setting
\[\theta= \sum_{\pi \in {{\bf P}}(m\gc)}e_{- \pi}H_{\pi}.\]
\noi Note that
$\theta(\gl) = \theta^\gl$.
The \v Sapovalov element in Theorem \ref{Shap} is constructed inductively using the next Lemma, see
for example  \cite{H2} Section 4.13 or \cite{M} Theorem 9.4.3.
\begin{lemma}\label{1768}
Let $\gc$ be a positive root, and $m$ a positive integer which is equal to 1 if $\gc$ is isotropic.
Suppose that  $\alpha \in
\Pi_{\rm even}, $ and set

$$\mu = s_\alpha\cdot \gl,\;\gc' = s_\alpha\gc,\;  p = (\mu + \rho, \alpha^\vee),\;q = (\gc, \alpha^\vee).$$
Assume that $p, q  \in \mathbb{N} \backslash \{0\}$ and
$\theta' \in U({\mathfrak n}^-)^{-m\gc'}$ is such that $v = \theta'v_\mu \in M(\mu)$
is a highest weight vector.
\noi Then there is a unique $\theta \in U({\mathfrak n}^-)^{-m\gc}$ such that

\begin{equation} \label{21nd}
e^{p + mq}_{- \alpha}\theta' = \theta e^p_{- \alpha}.
\end{equation}

\end{lemma}

%\bpf See for example  \cite{H2} Section 4.13 or \cite{M} Theorem 9.4.3.\epf
\noi Using Lemma \ref{1768} it follows that the coefficients of the \v Sapovalov elements $\gth_{\gc, m}$ are obtained by taking $\ttk$-linear combinations of
products of coefficients of the %\v Sapovalov elements
 $\gth_{\gc', m}$, with binomial  coefficients.

\section{Changing the Borel Subalgebra.}\label{sscbs}
Using adjacent Borel subalgebras (equivalently odd reflections), it is possible to give an alternative construction of \v Sapovalov elements corresponding to an isotropic root $\gc$, provided $\gc$ is a  simple root for some Borel subalgebra.  This condition always holds in type A, but for other types,  it is quite restrictive:  if $\fg = \osp(2m,2n+1)$, (an algebra that does not satisfy the restrictions imposed at the start of section 2) the assumption only holds
for roots of the form $\pm(\gep_i-\gd_j)$, while if $\fg = \osp(2m,2n)$ it holds only
for these roots and the root $\gep_m+\gd_n$.
\\ \\
Suppose that  %$\Pi$ is the  set of distinguished or anti-distinguished  set of simple roots $\Pi$. Suppose
$\fb$ is the distinguished  Borel subalgebra, and let $\mathfrak{b}'$
be another Borel subalgebra with the same even part as
$\mathfrak{b}.$ Consider  a  sequence of Borel subalgebras

\[\mathfrak{b} = \mathfrak{b}^{(0)}, \mathfrak{b}^{(1)}, \ldots,
 \mathfrak{b}^{(r)}. \]
  Assume there are isotropic roots $\ga_i$ such that $\fg^{\ga_i} \subset \fb^{(i-1)}, \fg^{-\ga_i} \subset \fb^{(i)}$  for $1 \leq i \leq r$, and  $\ga_1,\ldots,\ga_r$ are distinct positive roots of $\fb.$

\bt \label{rot} Set $F(\gc) = \{i|1 \leq i \leq r \mbox{ and } (\gc,\ga_i)=0  \}$. There is a nonzero $c\in \ttk$ such that for all $\gl \in \cH_\gc$,

\[e_{\ga_{1}}\ldots e_{\ga_r}e_{-\gc}e_{-\ga_{r}}\ldots e_{-\ga_1}v_\gl= c\prod_{i \in F(\gc)}(\gl+\gr,\ga_i)
\gth_\gc v_\gl.\]
 \et

\section{The Square of a \v Sapovalov Element.} \label{zzprod}
\noi When $\gamma$ is an isotropic root we write ${\cH}_{\gc}$ and $ \gth_{\gamma}$ in place of $\mathcal{H}_{\gamma,1}$ and $ \gth_{\gamma,1}$
 respectively. Here we record  an elementary but important property of the \v Sapovalov element $\theta_{\gamma}$ corresponding to such a root. %${\gamma}$.

\bt \label{zprod} If $\gl \in \cH_\gc,$ then $ \theta^2_{\gamma} v_\gl = 0.$ Equivalently,
  $\theta_{\gamma}(\gl-\gc)\theta_{\gamma}(\gl) =0.$\et
\noindent \bpf %It is enough to show this for all $\gl$ in a Zariski dense subset of $\cH_\gc$.
 We assume $\gc = w\gb$ for $\gb \in \Pi$ where $\gb$ is isotropic, and $w\in  W_{\rm even}$.
For a suitable Zariski dense subset $\gL$ of  ${\mathcal H}_{\beta}$, suppose $\gl \in
w\cdot \Lambda$.
The proof is by induction on $\ell(w).$ We can assume that
$w \neq 1.$ Suppose $w = s_\ga u$ with $\ell(u) = \ell(w) -1,$ and set $\gc' = s_\ga \gc.$ Replace $\mu$ with $\mu - \gc'$ and $\gl = s_\ga\cdot \mu$
with $s_\ga\cdot (\mu - \gc') = \gl - \gc$ in Equation (\ref{21nd}).
Then $p$ is replaced by $p + q$ and adopting the notation of \eqref{yam},

\[ e^{p + 2q}_{- \alpha}\theta^{\mu - \gc'}_{\gamma'} = \theta^{\lambda -\gc}_{\gamma}e^{p + q}_{- \alpha}.\]
 Combining this with Equation
(\ref{21nd}) and using induction we have

\[0 = e^{p+ 2q}_{-\ga}\theta^{\mu - \gc'}_{\gamma'}\theta^{\mu}_{\gamma'} =\theta^{\gl - \gc}_{\gamma} \theta^\gl_{\gamma}e^{p}_{-\ga} .\]
The
result follows since $e_{-\ga}$ is not a zero divisor in
$U({\mathfrak n}^-).$\epf

\section{Modules with Prescribed Characters.}\label{jaf}
We  introduce some new highest weight modules whose characters are given by generating functions for certain kinds of partitions.
If $X$ is a set of pairwise orthogonal isotropic positive roots,
set

\[{\bf P}_{X}(\eta) = \{\pi \in {\bf P}(\eta) | \pi(\alpha) = 0 \mbox{ for all } \ga \in X\}.\]
and ${\bf p}_{X}(\eta) =|{\bf P}_{X}(\eta)|$.
Set $p_X = \sum {\bf p}_X(\eta)\epsilon^{-\eta}$.
We have

\be \label{pfun}  p_X = \prod_{\alpha \in \Delta^+_{1}\backslash X} (1 +
\epsilon^{- \alpha})/ \prod_{\alpha \in \Delta^+_{0}} (1 -
\epsilon^{- \alpha}).\ee
If $X$ is empty, set $p = p_X$, and if $X =\{\ga\}$ is a singleton write

\be \label{sing}{\bf P}_{\ga}(\eta),\quad{\bf p}_{\ga}(\eta),\; \mbox{ and } \;p_\ga
\ee
instead of
${\bf P}_{X}(\eta),
{\bf p}_{X}(\eta),\mbox{ and } p_X.$
For a module $M$ in the BGG category $\cO$, the character of $M$ is defined by $\ch M= \sum_{\eta \in \fh^*}\dim_\ttk M^{\eta}\gep^\eta.$  Recall that the Verma module $M(\gl)$ has character $\gep^\gl p.$

\bt \label{newmod}
Suppose that  $X$ is an isotropic set of positive roots and $\lambda \in {\cH}_{\gc}$ for all $\gc\in X$. Then there
exists a factor module $M^{X}(\lambda)$ of
$M(\lambda)$ such that
\[ \ch M^{X}(\lambda) = \epsilon^{\lambda} p_{X}.\]
\et\noi
 If $X =\{\gc\}$ we write
 $M^{\gc}(\lambda)$ in place of $M^{X}(\lambda).$
\noi
The construction of the modules $M^X(\lambda)$ involves a process of deformation and specialization. First we  extend scalars to $A = \ttk[T]$ and $B = \ttk(T)$.  If $R$ is either of these algebras we set
$U(\fg)_R = U(\fg) \otimes R$. Choose $\xi \in \fh^*$ such that $(\xi, \gamma) = 0$ for all $\gc \in X$, and $(\xi, \ga^\vee) \notin \Z$ for all even roots $\ga$. Next consider the $U(\fg)_B$-module $M({{\widetilde{\lambda} }})_{B}$ with highest weight ${{\widetilde{\lambda} }} = \gl+T \xi$, and  form the factor module of $M({{\widetilde{\lambda} }})_{B}$ obtained
by setting $\gth_\gc v_{{{\widetilde{\lambda} }}}$ equal to zero for $\gc \in X$.
Then take a suitable $U(\fg)_{A}$-submodule of this factor module  and reduce mod $T$ to obtain the module $M^X(\lambda)$.
In more detail, we set

\[M^X({\widetilde{\lambda}})_B = M({\widetilde{\lambda}})_B/ \sum_{\gc \in X} U(\fg)_B \gth_\gc v_{{\widetilde{\lambda}}}.\]
Then $M^X({\widetilde{\lambda}})_B$ is  a $U(\fg)_B$-module generated by a highest weight vector
$v^X_{{\widetilde{\lambda}}}$
(the image of $v_{{\widetilde{\lambda}}}$)
with weight  $\widetilde{\lambda}$. Set
$M^X({\widetilde{\lambda}})_A = U(\fg)_Av^X_{{\widetilde{\lambda}}} \subset
M^X({\widetilde{\lambda}})_B,$ and
\be \label{yap}M^X({{\lambda}}) =M^X({\widetilde{\lambda}})_A /TM^X({\widetilde{\lambda}})_A.\ee
  Then
\[M^X({\widetilde{\lambda}})_A\ot_A B = M^X({\widetilde{\lambda}})_B.\]
Based on Corollary \ref{zoo} we can show
\bl \label{rats}Let $M=U(\fg)_Bv$ be a module with highest weight $\widetilde{\lambda}$ and highest weight vector $v$.  Suppose that $\gth_\gc v = 0$ for all $\gc \in X.$ Then\\ \\
\itema \;for all $\eta$ the weight space $M^{\widetilde{\lambda} -\eta}$ is spanned over $B$ by all $e_{-\pi}v$ where $\pi \in {\bf  P}_{X}(\eta)$.\\ \\
\itemb \;
$\dim_B M^{\widetilde{\lambda}- \eta} \leq {\bf p}_{X}(\eta)$.\el
\noi
{\it Proof of Theorem \ref{newmod} when $X =\{\gc\}$. } We set

%\be \label{vwer} v  = v_{{\widetilde{\lambda} }}\in {  M}_B = {M}({\widetilde{\lambda} } )_{B}, \quad u = \theta_{\gamma}v, \quad N_B = U(\mathfrak{g})_{B}u.\ee

\[ v  = v_{{\widetilde{\lambda} }}\in {  M}_B = {M}({\widetilde{\lambda} } )_{B}, \quad u = \theta_{\gamma}v, \quad N_B = U(\mathfrak{g})_{B}u.\]
For $R = A$ or $B$ we write $M^{\gc}({\widetilde{\lambda} })_R$ in place of $M^X({\widetilde{\lambda}})_R$.
Then   the module $M^{\gc}(\lambda)$ defined by \eqref{yap}
is generated by the image $\overline{v}$ of $v$ which is a highest weight
vector of weight $\lambda$.  Also
by Theorem \ref{zprod}, $\gth_\gc u = 0$, so we can apply Lemma \ref{rats} to both $N_B$ and $M_B/N_B$.
This gives

\be \label{dimmn}
\dim (M_B/N_B)^{\widetilde{\lambda}  - \eta} \le {\bf p}_{\gc}(\eta ),\quad
 \dim N_B^{\widetilde{\lambda}  - \eta}\le
 {\bf p}_{\gc}(\eta -\gc). \ee
Since

$${\bf p}(\eta) =  \dim {M_B}^{\widetilde{\lambda}  - \eta} = \dim (M_B/N_B)^{\widetilde{\lambda}  - \eta} + \dim N_B^{\widetilde{\lambda}  - \eta},$$
and ${\bf p}_{\gc}(\eta - \gc) +
 {\bf p}_{\gc}(\eta ) = {\bf p}(\eta),$ it follows that equality holds in (\ref{dimmn}).
Now we obtain the result from the following considerations applied to the weight spaces of the modules
$M^{\gc}({\widetilde{\lambda} })_{R}$ for $R = A, B$.   If $K$ is
an ${A}$-submodule of a finite dimensional $B$-module $L$ such that $K_{A} \otimes_{A} B =
L$, then $\dim_\ttk K/TK = \dim_B L.$
\hfill  $\Box$
%\\ \\
%We remark that
%\brs
\begin{rem} {\rm If $M'$ is the kernel of the natural map
$M(\lambda)\lra M^{\gamma}(\lambda)$, then
$U(\mathfrak{g})\theta_{\gamma}v_\gl\subseteq M',$ but  the inclusion can be strict. Indeed this happens when $\fg = \fsl(2,1)$ and $\gl =-\gr$ \cite{M17}, see also \cite{M} Exercise 10.5.4.}
\end{rem}
%\ers

\section{The Jantzen Sum Formula.} \label{s.1}
The Jantzen sum formula for a semisimple Lie algebra expresses the sum of the characters of the terms in the Jantzen filtration as a sum of characters  of Verma modules.  There is a version of the formula for contragredient Lie superalgebras in \cite{M} Theorem 10.3.1, but it contains some terms that are not characters of Verma modules. Here we see that these extra terms  are actually characters of the modules $M^{\gc}(\lambda)$ introduced in Theorem \ref{newmod}.
\\ \\
For $\lambda \in
\mathfrak{h}^*$ define

\[A(\lambda) =  \{ \alpha \in {\Delta}^+_{0} | (\lambda + \rho,
\alpha^\vee) \in \mathbb{N} \backslash \{0\} \}, \]
 \[B(\lambda) = \{ \alpha \in {\Delta} ^+_{1} | (\lambda + \rho,
\alpha) = 0 \} . \]
Now we state our improved version  of the Jantzen sum formula.
At the same time, rather than using characters as in \cite{M}, it is
useful to rewrite the result using the Grothendieck group $K(\cO)$ of the category $\cO$.  We define $K(\cO)$ to be the free abelian
group generated by the symbols $[L(\gl)]$ for $\gl \in \fh^*$.  If $M\in \cO$, the class of $M$ in $K(\cO)$ is defined as
$[M]=\sum_{\gl \in \fh^*} |{M}:L(\gl)|  [L(\lambda)] ,$
where $|{M}:L(\gl)| $ is the multiplicity of the composition factor $L(\gl)$ in $M$.

\bt \label{Jansum}
For all $\lambda \in \mathfrak{h}^*$
\be \label{lb} \sum_{i > 0} [{M}_{i}(\lambda)] = \sum_{\alpha \in A(\lambda)}[{M}(s_{\alpha}\cdot \lambda)] +
\sum_{\gc \in B(\lambda)} [M^{\gc}(\lambda -\gc)].\ee
\et

\bpf Combine Theorem \ref{newmod} with the result from \cite{M} Theorem 10.3.1.\epf
\noi The advantage of using this version of the formula is that $K(\cO)$ has a natural partial order.  For $A, B \in \cO$ we write $A
\ge B$ if $[A]-[B]$ is a linear combination of classes of simple modules with non-negative integer coefficients.  Clearly if $B$ is a subquotient of $A$
we have $[A] \ge [B]$.

\section{Orthogonal Isotropic Roots.}\label{ssuv.81}
We consider the structure of $M(\gl)$ when $B(\gl)$ consists of two orthogonal roots $\gc, \gc'$.  We say that $\gl\in \cH_{\gc} \cap \cH_{{\gc'}}$ is {\it weakly  generic} if $A(\lambda)= \emptyset,$ and
$B(\lambda)= \{\gc, {\gc'}\}$, and {\it generic} if $\mu$ is weakly generic for all
$\mu \in \gl +\mathbb{Z} \gc +\mathbb{Z} {\gc'}$. Choose $\xi \in \fh^*$ such that $(\xi, \gamma) = 0$ for all $\gc \in X$, and $(\xi, \ga) \neq 0$ for all even roots $\ga$.

\bl \label{fix} For $\gl\in \cH_{\gc} \cap \cH_{\gc'}$ there are only finitely many $c \in \ttk$ such that
\be \label{we1}\gth_\gc(\lambda +c\xi - \gc)\gth_{{\gc'}}(\gl+c\xi)=0\ee or
\be \label{we2}\gth_{{\gc'}}(\lambda +c\xi - \gc)\gth_\gc(\gl+c\xi)=0.\ee\el
\bpf Set $\widetilde{\lambda}= \gl +T\xi$.  It follows  from \eqref{hig} and Corollary \ref{zoo} that when $\gth_\gc(\widetilde{\lambda}- \gc)
\gth_{{\gc'}}(\widetilde{\lambda})v_{\widetilde{\lambda}}$
 is written as a $A$-linear combination of terms $e_{-\pi}v_{\widetilde{\lambda}}$, the coefficient of $e_{-\gc}e_{\gc'} v_{\widetilde{\lambda}}$ is a polynomial in $T$ of degree $d_\gc +d_{{\gc'}}$. Hence  there are only finitely many $c$ such that Equation (\ref{we1}) holds, and a similar argument applies to Equation (\ref{we2}).\epf
\noi It follows from Lemma \ref{fix} that the set

 \[\gL =\{\gl  \mbox{ generic in } \cH_{\gc} \cap \cH_{{\gc'}}|\gth_{\gc'}(\lambda - \gc)\gth_{\gc}(\gl)\neq 0\neq \gth_{\gc}(\lambda - {\gc'})\gth_{\gc'}(\gl)\}\]
is Zariski dense in $\cH_{{\gc'}} \cap \cH_{\gc}$.\\ \\
For $\gl \in \gL$, the Jantzen sum formula \eqref{lb} reads
\be \label{pig} \sum_{i > 0} [{M}_{i}(\gl)] =[M^{\gc}(\lambda -\gc)]  +[M^{{\gc'}}(\lambda -{\gc'})] .\ee
\noi Next we quote a result from \cite{M17} without proof.

\noi \bl \label{tri} For $\gl \in \gL$, $$|M(\gl):L(\gl-\gc-\gc')|=1.$$ \el
\noi Now based on equation \eqref{pig} we can show

\bl \label{hg}

\noi For $\gl \in \gL$ we have \\ \\
\itema \;If $|M(\lambda ):L(\mu)|>0,$ then $\mu \in \gl -\mathbb{N} \gc -\mathbb{N} \gc'$.\\ \\
\itemb \;$|M^{\gc}(\lambda):L(\gl)| = 1$.\\ \\
\itemc \;$|M^{{\gc'}}(\lambda-{\gc'}):L(\lambda -
r\gc)|=0,$ for $r >0.$\\ \\
\itemd \;$|M^{\gc}(\lambda):L(\lambda -
\gc)|=0.$\\ \\
\iteme \;$\sum_{i > 0} |{M}_{i}(\gl):L(\lambda-\gc)|=1.$\\ \\
\itemf \;$|M(\lambda):L(\lambda -r\gc)| = 0$ for $r\ge 2.$\\ \\
\itemg \;$|M(\gl):L(\gl-\gc-{\gc'})|=|{M}_{2}(\gl):L(\gl-\gc-{\gc'})|=1.$\\ \\
\itemh \;$|M^\gc(\gl):L(\gl-\gc-{\gc'})|=0.$\\ \\
\itemi \;$|M(\lambda):L(\gl-r\gc -s{\gc'})| =0$ if $r+s\ge 3.$
\el

\bpf (a) Follows from Equation (\ref{pig}) and induction on $\gl-\mu,$  while (b) holds since $M^{\gc}(\lambda)$ is an image of $M(\lambda).$
We note that (c) follows from (a).  To prove  (d) set $L= L(\gl-\gc)$, and suppose  for a contradiction, that $a =|M^{\gc}(\lambda):L|>0$. Then

\[a \le |M(\lambda):L|\le |M^{\gc}(\lambda-\gc):L| +|M^{{\gc'}}(\lambda-{\gc'}):L|=1,\]
by (b) and (c). Hence $a=|M(\lambda):L| =1,$  but then

$$|M(\lambda):L| = |M^{\gc}(\lambda):L|+|M^{\gc}(\lambda-\gc):L| = a+1,$$
since in the Grothendieck group $K(\cO)$, we have

\be \label{sop} [M(\lambda)] =  [M^{\gc}(\lambda)] + [M^{\gc}(\lambda-\gc)].\ee
Now (e) follows from (b), (c)  and (\ref{pig}), while (f) holds because for $r\ge 2,$

\[|M(\lambda):L(\lambda -r\gc)|\le |M^{\gc}(\lambda-\gc):L(\lambda -r\gc)| +|M^{\gc'}(\lambda-{\gc'}):L(\lambda -r\gc)|= 0,\]
by (c), (d) and induction.   Next set $L=L(\lambda -
\gc-{\gc'})$. To prove (g), it is enough, by Lemma \ref{tri} to show that $\sum_{i > 0} |{M}_{i}(\gl):L| = 2$.
First note that $|M^{\gc}(\lambda-2\gc):L| =0$ by (c), so $|M^{\gc}(\lambda-\gc):L|=  |M(\lambda-\gc):L|$ and likewise $|M^{\gc'}(\lambda-\gc'):L|=  |M(\lambda-\gc'):L|$.  So by (e),
$$\sum_{i > 0} |{M}_{i}(\gl):L| =  |M^{\gc}(\lambda-\gc):L| +|M^{{\gc'}}(\lambda-{\gc'}):L|=|M(\lambda-\gc):L| +|M(\lambda-\gc'):L|=2.$$
By  the proof of (g)
$|M^{\gc}(\lambda-\gc):L| = 1$, so (h) follows from (\ref{sop}).  Finally (i) follows from (\ref{pig}) and induction on $r+s$.  For example \[|M(\gl):L(\gl-\gc-2{\gc'})|\le
|M^\gc(\gl-\gc):L(\gl-\gc-2{\gc'})| + |M^{\gc'}(\gl-\gc'):L(\gl-\gc-2\gc')|=0,\]
using (f) and (g). \epf

\bt \label{hot} Suppose  $\gc$ and ${\gc'}$ are orthogonal isotropic roots, and $\gl\in \gL$.
% is   a generic point of $\cH_{\gc} \cap \cH_{{\gc'}}$.
Define

\[V_1  = U(\fg)\gth_\gc v_\gl,\quad  V_2 =U(\fg)\gth_{\gc'} v_\gl.
\]
Then the Jantzen filtration $\{M_i=M_i(\gl)\}_{i>0}$ on $M(\gl)$ is given by

\[ \label{lid1}{M}_{3} = 0, \quad {M}_{2} = V_1 \cap V_2 \cong L(\gl-\gc-{\gc'}), \quad {M}_{1}= V_1+ V_2.\]
Moreover $M_1/M_2$ is the direct sum of $L(\gl-\gc)$ and $L(\gl -\gc')$.
\et
\bpf By Lemma \ref{hg} the only composition factors of $M_1$ are $L(\gl-\gc), L(\gl-{\gc'})$, and $L(\gl-\gc-{\gc'})$ each with multiplicity one. By (e) and (g) in the Lemma $L(\gl-\gc-{\gc'})\subseteq {M}_{2}$ and
\[ \sum_{i > 0} [{M}_{i}] = [L(\lambda-\gc)]+[L(\lambda -\gc')]+2[L(\gl-\gc-{\gc'})].\]
 Note that $V_1$ and $V_2$ are generated by highest weight vectors of weights $\gl-\gc$ and $\gl-\gc'$.
The result follows since for $\gl \in \gL$,
$\gth_{\gc'}\gth_\gc v_\gl$ and $\gth_\gc\gth_ {\gc'}v_\gl$ are highest weight vectors with weight
$\gl-\gc-{\gc'}$.
\epf

\bc \label{kt1} There is a rational function $p$ of $\gl \in \cH_{\gc} \cap \cH_{{\gc'}}$ such that
\be \label{kt}\gth_{\gc'}(\lambda - \gc)\gth_\gc(\gl)=p(\gl)\gth_ \gc(\gl - {\gc'})  \gth_ {\gc'}(\gl).\ee
\ec

\bpf The highest weight vectors referred to in the last sentence of the proof are necessarily proportional.  Hence \eqref{kt} holds since the coefficients of  $e_{-\gc'}e_{-\gc}$ in these highest weight vectors (when written as linear combinations of the
$e_{-\pi}$ with  $\pi \in \bf{ P}(\gc+\gc')$) are polynomials in $\gl$.
\epf \noi
Corollary \ref{kt1} is applied in the proof of Theorem \ref{newmod} in the case $|X|\ge 2$.  We sketch the main new ingredient when $X = \{\gc,\gc'\}$ as in the Corollary.  In this case
$M({\widetilde{\lambda}})_B$ has a series of submodules

\be \label{cod} M({\widetilde{\lambda}})_B = W_0 \supset W_1  \supset W_2 \supset W_3 \supset W_4 =0,
\ee
where \[ W_1  = U(\fg)_B\gth_\gc v_{\widetilde{\lambda}} +U(\fg)_B\gth_{\gc'} v_{\widetilde{\lambda}} , \quad  W_2 =
U(\fg)_B\gth_\gc v_{\widetilde{\lambda}} , \quad   W_3 =
U(\fg)_B\gth_\gc \gth_{\gc'} v_{\widetilde{\lambda}}.
\]
From the Corollary and Theorem \ref{zprod} we deduce that
$\gth_\gc \gth_{\gc'} v_{\widetilde{\lambda}} \in W_2$,
$\gth_{\gc'} \gth_\gc v_{\widetilde{\lambda}} \in W_3$ and
$\gth_{\gc'} \gth_\gc \gth_{\gc'} v_{\widetilde{\lambda}}=0$. Thus from Theorem \ref{zprod} and Lemma \ref{rats}
we obtain bounds on the dimensions of the weight spaces of the factors in the series \eqref{cod} which are analogous to \eqref{dimmn}. The statement about characters follows as before.

\section{Non-Orthogonal Isotropic Roots.}
If $\gc,\gc'$ are  non-orthogonal isotropic roots, then
$\gc'=s_\ga \gc$ for some even root
$\ga$.  In this situation
we relate the \v Sapovalov elements for $\gc, \gc'$ and $\ga$. We set $\gth_{\ga,0} =1$ and $Q_0^+ =\sum_{\ga\in\Pi_{\rm even}}\N\ga$.
%Assume that $q =\pm 1$.
\bt  \label{man}Let  $\gc$ be a positive isotropic root and
$\alpha$ a non-isotropic root contained in ${Q}^+_{0}.$ Let  $v_{\gl}$ be a highest weight vector
in a Verma module with highest weight $\gl,$ and set $\gc' = s_\alpha\gc.$
Suppose $p = (\gl+\gr, \alpha^\vee) \in \mathbb{N}\backslash \{0\}$.
 Then

\itema \; If $(\gl+\gr,\gc') = 0$
we have
\be \label{pin}\theta_{\gc}\theta_{\ga,p} v_{\gl}= \theta_{\ga,p+1} \theta_{\gc'} v_{\gl}.\ee
\itemb \;
If $(\gl+\gr,\gc) = 0$, and $p-1\ge 0,$
we have

\be \label{pun}\theta_{\gc'}\theta_{\ga,p} v_{\gl}= \theta_{\ga,p-1} \theta_{\gc} v_{\gl}.\ee
\et
\bpf It suffices to prove (a) for all $\gl$ in the Zariski dense subset $\gL$ of
$\cH_{\gc'} \cap \cH_{{\ga,p}}$ given by
\[\gL=\{ \gl \in \cH_{\gc'} \cap \cH_{{\ga,p}}|A(\gl) = \{\ga\},\; B(\gl) = \{\gc'\}\}.\]
However for $\gl \in \gL$, $M(\gl)$ contains a unique highest weight vector of weight $s_\ga \cdot \gl -\gc.$ Thus, $\theta_{\gc}\theta_{\ga,p} v_{\gl}$ and $ \theta_{\ga,p+1} \theta_{\gc'} v_{\gl}$ are equal up to a  scalar multiple. If $\pi^0$ is the partition of $p\ga+\gc$ with $\pi^0(\gs)=0$ for all non-simple roots $\gs,$ it follows easily from the definition of \v Sapovalov elements, that $e_{-\pi^0}v_{\gl}$ occurs with coefficient equal to one in both $\theta_{\gc}\theta_{\ga,p} v_{\gl}$ and $ \theta_{\ga,p+1} \theta_{\gc'} v_{\gl}$,  and  this gives  the desired result.  The proof of (b) is similar.
\epf

\begin{rem}\rm{
The most interesting case of Equation (\ref{pin}) arises when $p=0$, since then we have an inclusion between submodules of a Verma module obtained by multiplying the highest weight vector $v_{\gl}$ by $\theta_{\gc}$ and $\theta_{\gc'}$.  Similarly the most interesting case of Equation (\ref{pun}) is when $p=1.$
}\er

\noi
{\it Acknowledgments.}\\
%\verb"\ack"
I thank Kevin  Coulembier for helpful comments.  The research that led to this paper
%The research of the author
was partly supported by  NSA Grant H98230-12-1-0249 and Simons Foundation grant 318264.
\\

\end{document}